\documentclass[12pt]{article}
\usepackage{amsthm,amstext,amsmath,amscd,amssymb,eucal,latexsym}

\newtheorem{theorem}{{\rm T\sc heorem}}[section]
  \newtheorem{lemma}[theorem]{{\rm L\sc emma}}
  \newtheorem{corollary}[theorem]{{\rm C\sc orollary}}
\newtheorem{proposition}[theorem]{{\rm P\sc roposition}}
\newtheorem{remark}[theorem]{{\rm R\sc emark}}

\newtheorem{definition}[theorem]{{\rm D\sc efinition}}

\def\Pn#1{{\bf P}^{#1}}
\def\Pnd#1{{\check{\bf P}}^{#1}}
\def\Gr #1#2{{\bf G}(#1,#2)} 
\def\LG #1#2{{\bf LG}(#1,#2)}
\def\PW{{\bf P}(W)}
\def\PWD{{\bf P}(W^{*})}

\begin{document}

\title{The abelian fibration on the Hilbert cube of a K3 surface of genus 9} 
\author{Atanas Iliev \ and \ Kristian Ranestad} 

\date{}
\maketitle

\footnotetext[0]{\ The 1-st author is partially supported by grant MM1106/2001
                  of the Bulgarian Foundation for Scientific Research}

\begin{abstract}
In this paper we construct an abelian fibration over ${\bf P}^3$ 
on the Hilbert cube of the primitive K3 surface of genus $9$. 
After the abelian fibration constructed by Mukai  on 
the Hilbert square on the primitive K3 surface of genus $5$, 
this is the second example where the abelian fibration on 
such ${\mathcal Hilb}_nS$ is directly constructed. 
Our example is also the first known abelian fibration on a 
Hilbert scheme ${\mathcal Hilb}_nS$ of a primitive
K3 surface $S$ which is not the Hilbert square of $S$;
the primitive $K3$ surfaces on the Hilbert square of which 
such a fibration exists are known by a recent result of 
Hassett and Tschinkel. 
\end{abstract}

\section{Introduction}

{\bf Generalities}. 
The smooth complex projective variety $X$ is a hyperk\"ahler manifold
if $X$ is simply connected and $H^0(\Omega_X^2) = {\bf C}\omega$ 
for an everywhere non-degenerate form $\omega$,
the symplectic form on $X$. In particular, by the non-degeneracy of $\omega$,
a hyperk\"ahler manifold is always even-dimensional and with a trivial canonical class, 
see e.g. \cite{Huy1} or \cite{Huy2} for a survey of the basic 
properties of hyperk\"ahler manifolds. 

A fibre space structure, or a {\it fibration} on the smooth projective 
manifold $X$ is a projective morphism, with connected equidimensional fibers,  
from $X$ onto a normal projective variety $Y$ such that $0 < dim(Y) < dim(X)$. 
By the theorems of Matsushita (see \cite{Mat}),
any fibration $f: X \rightarrow Y$ on a hyperk\"ahler $2n$-fold $X$ 
is always a Lagrangian abelian fibration, 
i.e. the general fiber $F_y = f^{-1}(y) \subset X$ 
must be an abelian $n$-fold which is also a Lagrangian 
submanifold of $X$ with respect to the form $\omega$;
in addition $Y$ has to be a Fano $n$-fold with the same 
Betti numbers as ${\bf P}^n$.  


The question when a hyperk\"ahler $2n$-fold $X$
admits a structure of a fibration in the above sense 
is among the basic problems for hyperk\"ahler manifolds 
(cf. \cite{Huy1}, p. 171).
An additional question is whether  
the base $Y$ of any such fibration is 
always the projective space ${\bf P}^n$, ibid. 
For a more detailed discussion on the problem, 
see the paper \cite{S} of J. Sawon. 

As shown by Beauville,    
the Hilbert schemes ${\mathcal Hilb}_n S$ 
of length $n$ zero-subschemes of smooth K3 surfaces $S$ are 
hyperk\"ahler manifolds, see \cite{Beau}.

In the special case when $S$ is an elliptic K3 surface,
the elliptic fibration $S \rightarrow {\bf P}^1$ on $S$ 
induces naturally a structure of an abelian fibration 
${\mathcal Hilb}_n S \rightarrow {\bf P}^n$ for any $n \ge 2$, 
see e.g. \cite{S}, Ex. 3.5.
Another candidate is Beauville's fibration 
${\mathcal Hilb}_g S \cdots> {\bf P}^g$ 
for a K3 surface $S$ containing a smooth curve of genus $g \ge 2$,
in particular if $S = S_{2g-2}$ is a primitive $K3$ surface of genus 
$g \ge 2$, ibid. Ex. 3.6. But Beauville's fibration is not regular,
i.e. it is not a fibration in the above sense. 

Until now, the only abelian fibration constructed directly 
on a Hilbert power of a primitive K3 surface $S = S_{2g-2}$,
$g \ge 2$ is the abelian fibration due to Mukai (\cite {Mu84})  on the Hilbert 
square of the primitive K3 surface $S_8$ of genus $5$. 
On the base of this example, and by using a deformation argument, 
Hassett and Tschinkel manage to prove
the existence of an abelian fibration over ${\bf P}^2$ on the Hilbert 
square of the primitive K3 surface of degree $2m^2$, for any $m \ge 2$, 
see \cite{HT}. 

\ 

{\bf The main result of the paper and a conjecture}. 
In this paper we construct an abelian fibration over ${\bf P}^3$ 
on the Hilbert cube of the general primitive K3 surface of genus $9$.
After Mukai's example, this is the second case 
where the abelian fibration on ${\mathcal Hilb}_n S$ 
is constructed directly. 
This is also the first known abelian fibration on a Hilbert 
power of a primitive $K3$ surface, which is not a Hilbert square.

In addition, we consider the conjecture that if $S$ is a primitive K3 surface 
of genus $g \ge 2$, then the Hilbert scheme $X = {\mathcal Hilb}_n S$ 
admits a structure of an abelian fibration if and only if 
$2g-2 = m^2(2n-2)$ for some integer $m \ge 2$.  This conjecture is posed 
as a question by several authors (cf. \cite{Huy1}, \cite{HT}), so we claim no originality.

The Hassett-Tschinkel result covers the cases $n = 2, m \ge 2$
of this conjecture, and Mukai's example corresponds 
to the particular case $n = m = 2$. Our example proves  
the conjecture in the case $n = 3, m = 2$. 

\

{\bf The construction of the abelian fibration in brief}. 
To explain our construction, we compare it with Mukai's 
construction of the abelian fibration
on the Hilbert square of the primitive K3 surface $S_8$ of genus $5$. 

The general $S = S_8$ is a complete intersection of three quadrics 
in ${\bf P}^5$.  Let ${\bf P}(H^0(I_S(2))) = \Pnd 2_S$ be the 
plane of quadratic equations of $S$ in ${\bf P}^5$. 
Inside $\Pnd 2_S$ the subset $\Delta_S$ of singular quadrics 
is a smooth plane sextic, and the double covering of $\Pnd 2_S$
branched along $\Delta_S$ defines uniquely a K3 surface $F_S$ of genus 
$2$, the dual K3 surface of $S$ (cf. \cite{Mu5}). 
Any subscheme $\xi \in {\mathcal Hilb}_2 S$ spans a line 
$l_{\xi} = <\xi> \subset {\bf P}^5$, and the set  
$L_{\xi} = {\bf P}(H^0(I_{S \cup l_{\xi}}(2))$ of quadrics 
containing $S$ and $l_{\xi}$ is a line in $\Pnd 2_S$. 
If ${\bf P}^2_S$ is the dual plane of $\Pnd 2_S$,
the association $f: \xi \mapsto  L_{\xi}$ defines a regular map 
$$
f: {\mathcal Hilb}_2 S \rightarrow {\bf P}^2_S. 
$$ 
In turn, a fixed line $L \in {\bf P}^2_S$ 
defines uniquely a 3-fold $X_L \supset S$ in  ${\bf P}^5$ 
as the common zero locus of the quadrics $Q \in L \subset \Pnd 2_S$.
By the definition of $f$, the $0$-schemes $\xi$ in the preimage 
$A_L = f^{-1}(L)$ are intersected on $S$ by the lines $L \subset X_L$,
thus the fiber $A_L$ is isomorphic to the Fano family 
$F(X_L)$ of lines on the threefold $X_L$. 
The theorems in \cite{Mat} imply that the map 
$f$ is an abelian fibration, in particular for the general $L$ 
the family $F(X_L)$ is an abelian surface. The last 
is classically known: the general $X_L$ is the quadratic complex 
of lines, and the family of lines $F(X_L)$ is the jacobian 
of a genus $2$ curve $F_L$, the dual curve to $X_L$;  
the last turns out to be a hyperplane section of $F_S$
-- the double covering $F_L$ of $L$ branched along the 6-tuple 
$\Delta_L = \Delta_S \cap L$, see e.g. \cite{D}.  
 
\medskip 

Let now $S = S_{16} \subset {\bf P}^9$ be a general primitive 
K3 surface of genus $9$,
and let ${\mathcal Hilb}_3 S$ be the Hilbert scheme of length 
$3$ zero-schemes $\xi \subset S$. 
Following the above example, we shall define step-by-step a map 
$$
f: {\mathcal Hilb}_3 S \rightarrow {\bf P}^3.
$$
The analog of the projective space ${\bf P}^5 \supset S_8$
is now the Lagrangian grassmannian $LG(3,6)$, the closed orbit 
of the projectivized irreducible $14$-dimensional representation $V_{14}$ 
of the symplectic group $Sp(3)$.
The Lagrangian grassmannian $\Sigma=\LG 36$ is a smooth $6$-fold 
of degree $16$ in ${\bf P}^{13} = {\bf P}(V_{14})$, 
and by the classification results of Mukai \cite{Mu7}, 
the general primitive K3 surface $S = S_{16}$ 
is a linear section of the Lagrangian grassmannian 
$LG(3,6) \subset {\bf P}^{13} = {\bf P}(V_{14})$ 
by a codimension $3$ subspace ${\bf P}^9$ identified with $Span(S)$. 

The projectivized dual representation $V_{14}^*$ of $Sp(3)$ 
has an open orbit in $\hat{\bf P}^{13} = {\bf P}(V_{14}^*)$, 
and the complement to this orbit is a quartic hypersurface 
$\hat{F}$. The space of linear equations 
$\hat{\bf P}^3_S = {\bf P}^{9,\perp}_S$ of $S \subset LG(3,6)$ 
is a projective $3$-space in $\hat{\bf P}^{13} = {\bf P}(V_{14}^*)$, 
which intersects $\hat{F}$ along a smooth quartic surface 
$F_S \subset \hat{\bf P}^3_S$, the $Sp(3)$-dual quartic of $S$, 
see \cite{IR}.  
The base of $f$ is now the dual $3$-space ${\bf P}^3_S = \hat{\bf P}^{3,*}_S$,
interpreted equivalently as the set of all $3$-folds 
$X_h = \Sigma \cap {\bf P}^{10}_h$, $h \in {\bf P}^3_S$ which lie 
inside $\Sigma = LG(3,6)$ and pass through $S$. 

\

{\bf Structure of the paper}. 
The crucial step in our construction is to find the analog of lines 
through 2 points on $S_8$. It turns out that these are twisted 
cubics on $\Sigma = LG(3,6)$, and the crucial property of these cubics is:
{\it Through the general triple $\xi$ of points on $\Sigma$ passes a unique 
twisted cubic $C_{\xi}$ that lies in $\Sigma$}.

\medskip

In Sections 2 and 3 we show that on the general K3 surface 
$S = S_{16} \subset LG(3,6)$, if $\xi$ 
is a $0$-scheme on $S$, then there exists on $LG(3,6)$ a unique  
connected rational cubic curve $C_{\xi}$ intersecting  
on $S$ the zero-scheme $\xi$. 
This identifies ${\mathcal Hilb}_3 S$ with the 6-fold Hilbert scheme
${\mathcal C}(S)$ 
of twisted cubic curves $C \subset LG(3,6)$ 
that lie in Fano $3$-folds $X_h, h \in {\bf P}^3_S$. 

If $F_S \subset \hat{\bf P}^3_S$ is the $Sp(3)$-dual 
quartic surface of $S$, 
any $h \in {\bf P}^3_S = \hat{\bf P}^{3,*}_S$ 
defines a hyperplane section $F_h \subset F_S$,
the dual plane quartic of $X_h$. 
For the general $h$, $F_h$ is a smooth plane quartic and $X_h$ 
is a smooth prime Fano $3$-fold of genus $9$; one can regard this $X_h$ 
as the analog of the general quadratic complex of lines $X_L$ 
through the K3 surface $S_8$.
For $h' \not= h''$ the families of cubics 
${\mathcal C}(X_{h'})$ and ${\mathcal C}(X_{h''})$  
do not intersect each other, thus giving a regular map 
$$
f: {\mathcal Hilb_3} S \cong {\mathcal C}(S) \rightarrow {\bf P}^3_S, 
$$
with $h = f(\xi)$ identified with the $10$-space 
${\bf P}^{10}_\xi = Span(S \cup C_{\xi})$. 
The results of \cite{Mat} imply that $f$ is an abelian fibration, 
and the construction of $f$ identifies the fiber $A_h = f^{-1}(h)$ 
of $f$ with the Hilbert scheme ${\mathcal C}(X_h)$ of twisted cubic curves on 
the $3$-fold $X_h$. 

\medskip

In Section 4 we show that for the general $h$ the family 
${\mathcal C}(X_h)$ is nothing else but the jacobian $J(F_h)$ 
of the dual plane quartic $F_h$ to the 3-fold $X_h$. 
Just as in the case for the quadratic complex of lines, 
the abelian threefold $A_h = {\mathcal C}(X_h)$ can 
be identified with the intermediate jacobian $J(X_h)$ 
of the Fano $3$-fold $X_h$. 

\medskip

At the end, 
in Section 5 we describe the group law on the general fiber 
$A_h$ of $f$, in the interpretation of $A_h$ as the 
Hilbert scheme ${\mathcal C}(X_h)$ of twisted cubic curves on $X_h$. 
More precisely, any twisted cubic $C_o \subset X_h$ 
defines on the abelian $3$-fold $A_h = {\mathcal C}(X_h)$ 
an additive group structure with $C_o = 0$,
and we identify which cubic curve on $X = X_h$ 
is the sum under this group operation of two general cubics $C', C'' \subset X$.
Notice that this is the analog of the Donagi's group 
law on the family $F(X_L)$ of lines on the 
$3$-fold quadratic complex of lines (see \cite{D}),  
identified with the general fiber of Mukai's 
abelian fibration on ${\mathcal Hilb}_2 S_8$. 
\vskip .5cm
We thank  Y. Prokhorov for the advice in computing
the multiplicity coefficients in Proposition \ref{bir},
and S. Popescu and D. Markushevich for the helpful 
discussions.

\section{Twisted cubic curves on $\LG 36$}

\
Let $V$ be $6$-dimensional vector space, let 
$\alpha\in\wedge^2V^{\ast}$ be a nondegenerate $2$-form on $V$.  
The natural linear map $d_{\alpha}:\wedge^3 V\to V$ induced by 
$\alpha$ has a $14$-dimensional kernel $W={\rm ker}d_{\alpha}$. 
Consider the Pl\"ucker embedding $\Gr 2V\subset {\bf P}(\wedge^3V)$. 
Then $\Sigma = \LG 3V={\bf P}(W)\cap \Gr 2V\subset {\bf P}(\wedge^{3}V)$ 
is the $6$-dimensional Grassmannian of 
Lagrangian planes in $\Pn 5={\bf P}(V)$ with respect to $\alpha$.  If 
$H\subset \Sigma$ is the hyperplane divisor in
this Pl\"ucker embedding, then the degree of $\Sigma$ is $H^6=16$, while the 
canonical divisor $K_{\Sigma}=-4H$. Thus a general linear section $\Pn 
{10}\cap \Sigma\subset \PW$ is a Fano threefold of genus $9$, while 
a general linear section $\Pn 
{9}\cap \Sigma\subset \PW$ is a $K3$ surface of genus $9$.  
\begin{proposition}\label{linsec} (Mukai \cite{Mu7}) A general 
prime Fano threefold of genus $9$ is the linear section $\Pn {10}\cap\LG 
3V$ for some $\Pn {10}\subset \PW$. 
A general $K3$-surface of genus $9$ is the linear section $\Pn {9}\cap\LG 
3V$ for some $\Pn {9}\subset \PW$.
\end{proposition}

Let $Sp(3)\subset SL(V)$ be the symplectic group of linear 
transformations that leaves $\alpha$ invariant.  Then $W$ is an 
irreducible $Sp(3)$-representation and the four orbit closures 
of $Sp(3)$ on $\PW$ are
$$
\Sigma \subset \Omega \subset F \subset \PW
$$
Their dimensions are $6,9,12$ and $13$.  In particular $F$ is a quartic 
hypersurface, the union of projective tangent spaces to $\Sigma$.
Similarly, there are four orbit closures of $Sp(3)$ in the dual space:
$$
\Sigma^{*} \subset \Omega^{*} \subset F^{*} \subset \PWD,
$$
isomorphic to the orbit closures in $\PW$.  In particular, the 
quartic hypersurface $F^{*}$ is the dual variety of $\Sigma$ and 
parametrizes tangent hyperplane sections to $\Sigma$.
If $X=\Sigma\cap L$ is a linear section of $\Sigma$, and
 $L^{\bot}\subset \PWD$ is the orthogonal linear space, then we denote 
 by $F_{X}$ the intersection
 $L^{\bot}\cap F^{*}$ and call it the $Sp(3)$-dual variety to $X$.

The four orbits in $\PW$ are characterized by secant properties of 
$\Sigma$ (cf. \cite{IR}):
\begin{proposition}\label{secants} Let $\omega \in \PW$,
then: 

\smallskip

(a) If $\omega \in \PW-\Omega$, then through $\omega$
passes a unique bisecant or tangent line $l_{\omega}$ to $\Sigma$. 
The line $l_{\omega}$ is tangent to $\Sigma$ if and only if 
$\omega \in F$. 

\smallskip

(b) If $\omega \in \Omega - \Sigma$ then the set of lines 
which pass through $\omega$ and are bisecant or tangent 
to $\Sigma$ sweep out a $4$-space 
$\Pn 4_{\omega} \subset \PW$,  
and the intersection $Q_{\omega} = \Pn 4_{\omega} \cap \Sigma$
is a smooth $3$-fold quadric. 
In the space $\Pn 5$, there exists a point $x = x(\omega)$ 
such that the quadric $Q_{\omega}=Q_{x(\omega)} \subset \Sigma$ 
coincides with the set of Lagrangian planes that 
pass through the point $x(\omega)$. 
\end{proposition}

Furthermore

\begin{lemma}\label{planes}

(a) A line not contained in $\Sigma$ intersect $\Sigma$
    in a $0$-scheme of length $\le 2$.

\smallskip
(b) $\Sigma$ does not contain planes, and a plane 
    that intersects $\Sigma$ along a conic section is contained in 
$\Omega$.   

\smallskip

(c) a plane $P\cong\Pn 2$ 
    such that the intersection scheme $P\cap \Sigma$
    contains a $0$-scheme $Z$ of length $3$ either intersects 
    $\Sigma$ exactly at $Z$ or $P \cap \Sigma$
    contains a line or a conic.

\smallskip

 (d) A $3$-space $P\cong\Pn 3$ 
     such that the intersection $P \cap \Sigma$ 
     contains a curve $C$ of degree $3$ 
     either intersects $\Sigma$ exactly along $C$ and $C$ is defined by a 
     determinantal net of quadrics or $P\subset \Omega$ and 
     $P\cap \Sigma$ is a quadric surface
     of rank $\ge 3$. 
     
     \end{lemma}
     
     \begin{proof}

Since $\Sigma$ is defined by quadrics, (a) 
     follows.  
The planes that are parameterized by a conic have a common point 
$p\in\Pn 5$, 
so the conic lies in the $3$-dimensional smooth quadric  $Q_p$ 
(\ref{secants}).  Since this quadric is smooth, it and therefore also $\Sigma$, 
contains no planes, and the conic is a plane section of $Q_{p}$.  Hence 
$P\subset <Q_p>\subset\Omega$ and (b) follows.
     For (c), consider a plane $P$ whose intersection with $\Sigma$ 
contains 
     a scheme $Z$ of length $4$, and assume that the intersection is 
     zero-dimensional.  Then $Z$ must be a complete intersection of 
two conics, and the intersection $P\cap \Sigma$ is precisely $Z$.
     Let $S$ be a general $\Pn 9$ that 
     contains $P$, and let 
     $S = \Sigma \cap \Pn 9$.  Then, by Bertini, $S$ is 
     a smooth surface. In fact $S$ is a $K3$-surface.   Since $Z$ 
     define dependent conditions on hyperplanes $h$, there is a 
nontrivial 
     extension of $I_{Z,S}$ by $O_{S}$ which define a rank $2$ sheaf 
     $E$ on $S$.  Since no length three subscheme on $S$ is contained 
     in a line, $E$ is a vector bundle (cf. \cite{OSS}), and since no plane 
     intersects $S$ in a scheme of length five $E$ is base point 
     free.  Therefore a general section of $E$ is a smooth subscheme 
     of length four contained in a plane.  But through this subscheme 
     there are two lines that meet in a point outside $S$, 
     contradicting the above proposition.

    Let $P$ be a $\Pn 3$ that intersects 
     $\Sigma$ in a curve $C$. By (c) this curve has degree at most 
$3$, so 
     for (d) we may assume that the degree is three.  Again by (c), 
     the intersection is pure: There are no zero-dimensional 
     components. 
     Now, $C$ cannot be a plane curve by part a).   Furthermore, since $\Sigma $ is 
     the intersection of quadrics, $C$ is contained in at least three 
     independent quadrics.  Pick two without common component, then 
for degree reasons alone 
     they link $C$ to a line and $C$ is defined by a 
     determinantal net of quadrics.  If $P$ intersects $\Sigma$ in a 
surface, this surface, by (a) is a quadric of rank at least $3$, and 
by (b) lies in $\Omega$.
      
   \end{proof}

\begin{lemma}\label{2x3} Let $M$ be a $2\times 3$-matrix of linear forms in $\Pn 
3$ whose rank $1$ locus is a curve $C_M$ of degree $3$, and let $p\in 
\Pn 3\setminus C_M$.  Then there is a unique line passing through 
$p$ that intersects $C_M$ in a scheme of length $2$, unless $p$ lies 
in a plane that intersects $C_M$ in a curve of degree $2$.
\end{lemma}

\begin{proof}  Let $M(p)$ be $M$ evaluated in $p$.   Then $M(p)$ has 
rank $2$, so let $(a_1,a_2,a_3)$ be the unique (up to scalar) 
solution to $M(p)\cdot a=0$.  Then $M\cdot a$ defines two linear 
forms that vanish on $p$. If they are independent, they define a line 
that intersect $C_M$ in a scheme of length $2$, and if they are 
dependent, then they define a plane that intersect $C_M$ in a curve 
of degree $2$. The uniqueness of the line in the first case follows by 
construction. \end{proof}
\

\begin{lemma}\label{lgsecants}  Let $\Pn 9\subset \PW$ be 
general.  In particular assume that $S = \Sigma \cap \Pn 9$ 
is a smooth $K3$ surface of genus $9$ with no rational curve of degree less than $4$.
    Then:      
     
\smallskip

(a) A line $l$ can intersect $S$ in at most a $0$-scheme 
    of length $\le 2$.

\smallskip
     
(b) If a plane $\Pn 2$ intersects $S$ 
     in a scheme containing a $0$-scheme $\xi$ of
     length $3$ then $\Pn 2 = <\xi >$ 
     and $\Pn 2 \cap S = \xi$. 

\smallskip

(c) If a $3$-space $\Pn 3$ is such that the 
    intersection scheme $\Pn 3 \cap S$ contains 
    a $0$-scheme $\xi$ of length $3$ and the intersection 
    $\Pn 3 \cap \Sigma$ 
    contains a curve $C$ of degree $3$, 
    then $\Pn 3 \cap \Pn 9 = <\xi >$,  
    $\Pn 3 \cap S = \xi$ and $\Pn 3 \cap \Sigma = C$;
    here $\Pn 9 = <S>$.   
\end{lemma}
\
        
\begin{proof}  
 
(a) follows immediately from Lemma \ref{planes}.
By (a), $<\xi > = \Pn 2$. If the intersection scheme $Z = \Pn 2 \cap S \supset 
\xi$
   contains $\xi$ properly, then by Lemma \ref{planes}
    $Z$ will contain a line or a conic.  
    But $S$ contains no curves of degree less than $4$, so (b) follows.   
      
(c) Since $\Pn 3 \cap S \supset \xi$ then $\Pn 3$ contains 
the 
    plane $\Pn 2_{\xi} = <\xi > \subset \Pn 9$, by (b). 
    Therefore either $\Pn 3 \subset \Pn 9 = <S>$ 
    or $\Pn 3 \cap \Pn 9 = <\xi >$. 
    But if $\Pn 3 \subset \Pn 9$ then the curve
    $C \subset \Pn 3 \cap \Sigma$ will be contained in 
    $S = \Sigma \cap \Pn 9$, which is impossible since 
    by assumption $S = S_{16}$ does not contain curves. of degree $3$. 
    Therefore $\Pn 3 \cap \Pn 9 = <\xi >$ and 
    $\Pn 3 \cap S = <\xi > \cap S = \xi$.   
    
    If $\Pn 3 \cap \Sigma$ contains more than the curve $C$ 
    then by Lemma \ref{planes}, the intersection $Q = \Pn 3 \cap \Sigma$ 
    will be a quadric surface of rank at least $3$. 
    But then $S\cap Q$ is a conic, contradicting the assumption that 
    $S$ has no curve of degree less than $4$. 
     \end{proof}

\
\begin{definition} Let ${\mathcal Hilb}_{3t+1}(\Sigma)$ be the Hilbert 
scheme of twisted 
cubic curves contained in $\Sigma$ and let ${\mathcal 
Hilb}_{3}(\Sigma)$ be the Hilbert 
scheme of length three subschemes of $\Sigma$. \end{definition}

Since $\Sigma$ contains no planes, every member $C$ of 
${\mathcal Hilb}_{3t+1}(\Sigma)$ is a curve 
of degree $3$ defined by a determinantal net of quadrics:  
$C$ lies in at least 3 quadrics, and since there are no planes in 
$\Sigma$ two general ones intersect in a curve of degree $4$ that 
links $C$ to a line.

Consider the incidence
$$I_{3}=\{(\xi, C)|\xi\subset C\}\subset{\mathcal Hilb}_{3}(\Sigma)\times {\mathcal 
H}_{3t+1}(\Sigma) $$
and the restriction to $S$:

$$I_{3}(S)=\{(\xi, C)\in I_{3}|\xi\subset S\}$$

\begin{definition} Let ${\mathcal C}(S)\subset {\mathcal 
Hilb}_{3t+1}(\Sigma)$ be the image of the projection $$I_{3}(S)\to {\mathcal 
Hilb}_{3t+1}(\Sigma),$$ i.e.
the Hilbert scheme of twisted cubic curves in $\Sigma$ 
that intersect $S$ in a scheme of length three.\end{definition}

\


\begin{corollary}\label{CS}
    Let $S=\Pn 9\cap \Sigma\in\PW$ be a smooth linear section with no 
    rational curves of degree less than four.

(a) ${\mathcal C}(S)\subset{\mathcal 
Hilb}_{3t+1}(\Sigma)$ is a closed subscheme.

\smallskip

(b)  The intersection map 
$$
\sigma: {\mathcal C}(S) \rightarrow {\mathcal Hilb}_3 S, 
\ C \mapsto C \cap S 
$$
is well defined on any $C \in {\mathcal C}(S)$.  

\end{corollary}
\begin{proof}  Consider the map ${\mathcal 
Hilb}_{3t+1}(\Sigma)\to \Gr 4{14}$ defined by $C\mapsto <C>$.  By 
Lemma \ref{lgsecants}, the subset ${\mathcal C}(S)$ is simply the 
pullback under this map of the closed variety of spaces that intersect $\Pn 9$ in 
codimension one.  Therefore (a) follows, while (b) follows 
since no component of $C$ is contained in $S$. 
\end{proof}
\

\begin{lemma}\label{sur}
The intersection map $\sigma: {\mathcal C}(S) \rightarrow {\mathcal Hilb}_3 S$ 
is surjective. 
\end{lemma}

\begin{proof}
By Corollary \ref{CS} the image $\sigma({\mathcal C}(S))$ is a closed 
subscheme of $Hilb_3 S$. 
Since $Hilb_3 S$ is irreducible, to prove that $\sigma$ is surjective 
it is enough to see that $\sigma$ is dominant. 

For a point $x \in S$ denote by $P^2_x \subset \Pn 5={\bf P}(V)$ 
the Lagrangian plane of $x$, and let 

\medskip
 
\centerline{${\mathcal U} = \{ \xi \in {\mathcal Hilb}_3 S: \xi = x + y + z$ 
is reduced and such that} 

\centerline{\ \ \ \ \ \ \ \ \ \ \ \ \ \ \ \ \ \ the Lagrangian planes 
$P^2_x$, $P^2_y$ and  $P^2_z$ are mutually disjoint $\}$.} 

\medskip

\noindent
Clearly ${\mathcal U} \subset {\mathcal Hilb}_3 S$ is open and dense,
so it rests to see that 

\medskip 

{\it For any $\xi = x+y+z \in {\mathcal U}$ 
there exists a smooth twisted cubic $C \subset \Sigma$ 
that passes through $x$, $y$ and $z$}. 

\medskip

Let $U_o$, $U_{\infty}$
and $U_y$ be the Lagrangian $3$-spaces of $x,y$ and $z$ in 
the $6$-space $V$,
i.e. $P^2_x = {\bf P}(U_0)$, 
$\Pn 2_z = {\bf P}(U_{\infty})$
and $\Pn 2_y = {\bf P}(U)$.

Since $\Pn 2_x$ and $\Pn 2_z$ do not intersect each other,
we may write:
$$
V = U_0 \oplus U_{\infty}.
$$
Furthermore, since $U_0$ and $U_{\infty}$ are Lagrangian,
we may choose coordinates $(e_i,x_i)$ on $U_0$
and $(e_{3+i},x_{3+i})$ on $U_{\infty}$,
$i = 1,2,3$
such that in these coordinates the form
$\alpha$ can be written as
$$
\alpha = x_{1}\wedge x_{4}+x_{2}\wedge x_{5}+x_{3}\wedge x_{6}.
$$
Since $\Pn 2_y = {\bf P}(U)$ is also Lagrangian,
and since $\Pn 2_y$ does not intersect $\Pn 2_z$,
there exists a symmetric non-singular $3 \times 3$ matrix $B$ 
such that the Pl\"ucker coordinates of $y$ in the
system $(e_i,x_i)$, $i = 1,...,6$ are uniquely
written in the form
$$
y = exp(B) = (1:B:\wedge^2 B:det(B)).
$$
  For a parameter $t\in {\bf 
C}$, the 
matrix $tB$ correspond to a Lagrangian plane that does not intersect 
$\Pn 2_{z}$.
Thus 
$$
C = C_{x,y,z} = 
\ \{ exp(t) = (1: tB:t^2\wedge^2 B: t^3\det(B)), \
t \in {\bf C} \cup {\infty} \}
$$
is a smooth twisted cubic on $\Sigma$ through 
$x = exp(0)$, $y = exp(1)$ and  $z = exp(\infty)$.
\end{proof}

\

\begin{proposition}\label{iso} 
Let $S=\Pn 9\cap \Sigma$ be a smooth linear section with no 
    rational curves of degree less than four, and let ${\mathcal C}(S) $ be 
    the Hilbert scheme of twisted cubic curves on $\Sigma$ 
    that intersect $S$ in a scheme of length three.
    Then the restriction map
$$\sigma: {\mathcal C}(S) \rightarrow {\mathcal Hilb}_3 S\quad C\mapsto C\cap S$$ 
is an isomorphism. In particular ${\mathcal C}(S)$ is a smooth projective variety. 
\end{proposition}
\begin{proof}
First we show that $\sigma$ is bijective.
Let $\xi\in {\mathcal Hilb}_3(S)$ let $P$ be the span of $\xi$ and assume that 
$C_1,C_2\in{\mathcal C}(S)$ with $C_i\cap S=\xi$ for $i=1,2$.

 Let $q\in P$ be a general point. Since $S$ contains no plane curves, 
$P$ is not contained in $\Omega$, so we may assume that $q$ does not 
lie in $\Omega$.  Therefore there is a unique line through $p$ that 
intersects $\Sigma$ in a scheme of length $2$. On the other hand for 
each $i$, there is a unique line through $q$ that meets $C_i$ in a 
scheme of length $2$.  Hence these lines must coincide, and lie in 
$P$.  But the only way the general point $q$ in $P$ can lie on  a 
unique line through $\xi$, that intersects $\xi$  in a scheme of 
length $2$, is that $\xi$ is the first order neighborhood of a point, 
i.e. when $P$ is the tangent plane to $S$ at the support $p$ of $\xi$.
In this case both $C_1$ and $C_2$ must be singular at $p$ and have at 
least one line component.  In fact, if the tangent cone to $C_i$ is 
planar, then this plane must contain the line component through $p$.  
Since $S$ does not contain lines, the tangent cone must span the 
$3$-space, and each $C_i$ is contained in the tangent cone to 
$\Sigma$ at $p$.  But this tangent cone is the cone over a Veronese 
surface (cf. \cite{IR}). In the projection from $p$, the plane $P$ is mapped to a 
line $L$ in the space of the Veronese surface.  The line $L$ does not 
intersect the Veronese surface, so it is not contained in the cubic 
hypersurface secant variety of the Veronese surface. It is a 
well known classical fact that there is a unique plane through $L$ 
that meet $V$ in a scheme of length $3$.  This plane and the plane 
$P$ spans a $\Pn 3$ that intersects $\Sigma$ along three lines 
through $p$.  Therefore, also in this case  $C_1$ and $C_2$ must 
coincide.

Next,  we consider ramification of the map $\sigma$.  
Consider the morphism $g:{\mathcal C}(S)\to \Gr 4{14}$, defined by 
$[C]\mapsto [<C>]$,
and similarly $g_{3}:{\mathcal Hilb}_3 S\to \Gr 3{10}$ defined by 
$[\xi]\mapsto [<\xi>]$. By Lemma \ref{lgsecants}, 
both $g$ and $g_{3}$ are embeddings.
Therefore the restriction map $\sigma$ factors through $g$, the restriction 
map $[\Pn 3]\mapsto [\Pn 3\cap\Pn 9]$ 
and the inverse of $g^{-1}_{3}$ restricted to Im$g_{3}$.  
Hence a point of ramification for $\sigma$ is also a 
point of ramification for the restriction map $[<C>]=[\Pn 3]\mapsto [\Pn 
3\cap\Pn 9]=[<\xi>]$.  
If $[<C>]$ is such a point, then there is a tangent line to Im$g$ at 
the point $[<C>]\in \Gr 4{14}$ that is collapsed by the restriction map.  
But this occurs only if the tangent line is contained in $\Gr 4{14}$ 
and parameterizes $\Pn 3$s through $<\xi>$ in a $\Pn 4$.  
Since $<C>\cap \Sigma =C$, and $\sigma$ is a 
bijective morphism, no other $\Pn 3$ in this pencil can intersect $\Sigma$ 
in a twisted cubic curve. 
On the other hand, the ramification means that the doubling of $<C>$ 
in $\Pn 4$ intersect $\Sigma$ in a doubling $D$ of $C$,
i.e. $D$ is a nonreduced curve of degree $6$.  
For each point $p\in C$ consider the span $P_{p}$ of the tangent cone to $D$ at $p$.  
On the one hand, any line in $P_{p}$ through $p$ is tangent to $D$. 
On the other hand, the dimension of $P_{p}$ is one more than 
the dimension of the span of the tangent cone to $C$ at $p$.
Notice that if $p$ and $q$ are distinct points on $C$, 
and the line  $\overline{pq}$ is not a component of $C$, 
then $P_{p}\cap P_{q}\cap \overline{pq}=\emptyset$; 
otherwise $\overline{pq}$ is tangent 
to $D$ at $p$ or $q$ and is therefore contained in $C$ by Lemma \ref{planes}. 
  

Now, either $C$ is reduced, and two general tangent lines span $<C>$, 
or $C$ has a singular point whose tangent cone spans $<C>$.

 In the first case, let $p$ and $q$ are smooth points on $C$ such 
 that their tangent lines to $C$ span $\Pn 3$, then $P_{p}$ and $P_{q}$  
 are planes and $P_{p}\cap P_{q}$ meet in a point $r$.  
 The lines $\overline{pr}$ and  $\overline{qr}$ are both tangent to $D$  
 so the plane spanned by 
$p, q, r$ is by Lemma \ref{planes} contained in $\Omega$. 
 As $p$ and $q$ moves, the lines $\overline{pq}$ fill $<C>$, 
 so $<C>\subset\Omega$ against the assumption.
 
In the second case, let $p$ be a singular point on $C$ such that the tangent 
cone to $C$ at $p$ spans $<C>$, then $P_{p}=\Pn 4$ is contained in the 
tangent cone to $\Sigma$ at $p$. In particular $D$ is contained in 
this tangent cone.  But this tangent cone is a cone over a Veronese 
surface, so if $P_{p}$ intersects this cone in a curve of degree $6$, 
the intersection cannot be proper. In fact, the only improper 
intersections of a $\Pn 4$ with a cone over a Veronese surface that 
contains the vertex, is the union of a quadric cone surface and a 
line.  But in our case, the quadric cone would then have to be 
contained in $<C>$, contradicting the assumption.
 
 This concludes the argument that $\sigma$ is bijective and unramified, 
 i.e. an isomorphism.

\end{proof}

\
\section{The abelian fibration on ${\mathcal Hilb}_3 S$}

The Hilbert scheme ${\mathcal C}(S)$ of twisted cubic curves on $\Sigma$
 that intersect $S$ in a scheme of length three, admits a natural 
 fibration over $\Pn 3$.  
 Composed with the inverse of the restriction map 
 $\sigma: {\mathcal C}(S)\to {\mathcal Hilb}_3 S$, 
 we obtain an fibration on  ${\mathcal Hilb}_3 S$.  
 In this section we describe this fibration and explain how it may also 
 be obtained by a vector bundle approach as suggested in the introduction.  
 An upshot of the alternative approach is the
 existence of genus three curves in the fibers (cf. Proposition \ref{ajm}).
 \
  Let $C\subset \Sigma$ be a twisted cubic curve that intersect $S$ 
  in a scheme of length three.  
  Then the union $C\cup S$ spans a $\Pn {10}(C)\subset <\Sigma>=\PW$.  
  Since $$\{\Pn {10}|S\subset {\bf 
P}^{10}\subset\PW\}\cong {\bf P}_S^{3}$$
the assignment $ C\mapsto \Pn {10}(C)$
defines a map
$$\tilde{f}:{\mathcal C}(S)\to {\bf P}_S^3.$$
For any $h\in\Pn 3_{S}$, let  $X_h = \Sigma\cap\Pn {10}_{h} \supset S$. 
Then $X_{h}$ is a threefold, in fact, the general $X_h \supset S$ is a 
smooth Fano threefold of genus $9$. 
The fiber  $\tilde{f}^{-1}(h)$ of the map  $\tilde{f}$ is 
$$
{\mathcal C}(X_{h})=\{C\in {\mathcal Hilb}_{3t+1}(X_{h})\}.
$$

Assume that $S$ contains no curves of degree $\le 3$.  Then the map $\tilde{f}$ is 
a morphism.  On the other hand the restriction map 
$\sigma: {\mathcal C}(S)\to {\mathcal Hilb}_3 S$ is in his case an isomorphism, so
the inverse of $\sigma$ composes with $\tilde{f}$ to define a morphism
$$f=\tilde{f}\circ\sigma^{-1}:{\mathcal Hilb}_3 S\to  {\bf P}_S^3.$$

By the theorems of Matsushita (cf. \cite{Mat}), 
$f$ is a fibration on the hyperk\"ahler $6$-fold $Hilb_3 S$, 
in particular $f$ is a Lagrangian {\it abelian fibration}, 
i.e. the general fiber $A_h$ is an abelian $3$-fold
which is a Lagrangian submanifold with respect to the nondegenerate 
$2$-form on the hyperk\"ahler $6$-fold ${\mathcal Hilb}_3 S$.  
Thus 

\

\begin{theorem}\label{Af}
    Let $S=\Pn 9\cap \LG 3V\subset \PW$ be a smooth linear section with no 
    rational curves of degree less than four, then
${\mathcal Hilb}_3 S$ admits a fibration 
$f : {\mathcal Hilb}_3 S \rightarrow \Pn 3_{S}$.

\smallskip 

(a) For any $h \in \Pn 3_{S}$ 
the fiber $f^{-1}(h) = A_h=\sigma({\mathcal C}(X_h)$ 
where ${\mathcal C}(X_h) $ is the Hilbert scheme of twisted cubic curves 
on the threefold $X_h = \Sigma\cap\Pn {10}_{h} \supset S$ and $\sigma$ 
is the isomorphism above, defined by
$C\mapsto C\cap \Pn 9$. 

\smallskip

(b) For the general $h \in \Pn 3_{S}$ 
the fiber $A_h$ of $f$ is an abelian $3$-fold
which is a Lagrangian submanifold of the hyperk\"ahler 
$6$-fold ${\mathcal Hilb}_3 S$,
i.e. 
$$f : {\mathcal Hilb}_3 S \rightarrow \Pn 3_{S}$$
is a Lagrangian \underline{abelian fibration} on ${\mathcal Hilb}_3 S$.  
 \end{theorem}

By the general choice of $S$, the general $X_h \supset S$ 
is a general smooth prime Fano $3$-fold of genus $9$ (cf. Proposition 
\ref{linsec}), 
so the above theorem implies

\begin{corollary}\label{cxabel}
The Hilbert scheme of twisted cubic curves 
${\mathcal C}(X)={\mathcal Hilb}_{3t+1}(X)$ of the general prime Fano threefold 
$X$ of genus $9$ is isomorphic to a $3$-dimensional abelian variety.   
\end{corollary}

In the next section we shall show that the abelian variety ${\mathcal C}(X_h) \cong A_h$
 coincides with the Jacobian of the $Sp(3)$-dual plane quartic curve $F_{X}$ to $X$.   

\

Consider also an alternative approach.

First recall from \cite{IR} section 3.2 that zeros of a 
general $2$-form $\beta\in \wedge^2V^*$ 
on $\Sigma=LG(3,V)$ is a Segre threefold $\Pn 1\times\Pn 1\times\Pn 1$.  
The two $2$-forms $\alpha$ and $\beta$ define a unique triple of conjugate 
lines with respect to the two forms.
More precisely, $\alpha$ defines a correlation $L=L_{\alpha}$ which 
assigns to a point a hyperplane in $\Pn 5 $ through the point.  
This correlation is symmetric, 
in the sense that $y\in L(x)$ if and only if $x\in L(y)$.
 A triple of lines $L_{1}, L_{2}, L_{3}$ in $\Pn 5$
are conjugate with respect to $\alpha$ 
if $\cap_{x\in L_i}L(x)=<L_j,L_k>$ when ${1,2,3}={i,j,k}$. 
 If $\alpha$ and $\beta$ are nondegenerate, and $L_{1}, L_{2}, L_{3}$ 
 are conjugate w.r.t. 
 both $\alpha$ and $\beta$, then the planes in $\Pn 5$ that are Lagrangian 
with respect to both forms are precisely the planes that meet the three lines.   
Clearly these planes are parameterized by a Segre threefold in $\Sigma$, 
and any Segre threefold appear this way.

\begin{lemma}\label{unisegre}  A general twisted cubic curve $C$ in $\Sigma$  
and a general tangent hyperplane section $H_t$ to $\Sigma$
 that contains $C$, determines uniquely a Segre threefold 
 $\Pn 1\times\Pn 1\times\Pn 1$ in $\Sigma$ that contains $C$ and this 
 threefold lies in $H_t$.

The first general means that $C$ parameterizes planes of a $3$-fold 
$V_C\cong \Pn 1\times \Pn 2 \subset \Pn 5$. 
The second general means that the plane $P$ of the point of tangency 
of $H_t$ meets this $3$-fold in $3$ distinct points.
\end{lemma}

\begin{proof}  The $3$ points $x_1, x_2, x_3$ in $P\cap V_C$ 
determines three lines $L_1, L_2, L_3$ in 
$V_C$ that meet each plane in $V_C$ in a point.  
It is enough to prove that these three lines are conjugate 
w.r.t. $\alpha$ since uniqueness follows from the above.  
Notice that for each point $x$ in $V_C$ 
the hyperplane $L(x)$ intersects $V_C$ in the plane 
$P_x$ (through $x$) and a quadric surface $Q_x=\Pn 1\times \Pn 1$.  
So for each point in $Q_x$ 
there is a unique line in $Q_x$ that meet every plane in $V_C$ in a point.
Consider the three points $x_i$ and the corresponding quadrics $Q_i=\Pn 1\times \Pn 1$.
Since $x_i \in P$,  the hyperplane $L(x_i)$ contains $P$, so  $x_j$ and $x_k$ 
also lies in $L(x_i)$.  Therefore $x_j,x_k\in Q_i$, and $L_j, L_k$ lies in $Q_i$.
Let $x\in L(x_j)$ be a general point, then $L(x)$ by symmetry contains $x_i$, 
so $L_i$ must lie in the quadric $Q_x$.  
In particular  the hyperplane $L(x)$
 contains $L_i$. Therefore the three lines $L_1, L_2, L_3$ are 
 conjugate with respect to $\alpha$.
The corresponding Segre threefold in $\Sigma$ contains $C$ and the point of $P$, 
so by [IR] Lemma 3.2.5 it is contained in $H_t$.  
\end{proof}

Recall from \cite{IR} the following proposition 
\begin{proposition}\label{vb}
For each $t\in F^*\setminus 
\Omega^*\subset\PWD$, the tangent hyperplane section $H_t\subset \Sigma$ projects 
from its point $v(t)$ of tangency to a linear section of $G(2,6)$, 
hence it admits a rank $2$ vector bundle $E_t$.  
Every Segre threefolds in $\Sigma$ that pass through $v(t)$ is the zero locus
 of a section of this vector bundle, and the general section of $E_t$ is of this kind.  
\end{proposition}

\begin{remark}{\rm For the general $K3$-surface section $S$ of $\Sigma$, the set of 
tangent hyperplane sections $H_t$ that contain $S$ form the $Sp(3)$-dual quartic 
surface $T=F_{S}$, hence the same surface $T$ parameterizes rank $2$ vector 
bundles $E$ on $S$ with determinant $H$ and $H^0(S,E)=6$.  The general 
section of a vector bundle $E$ vanishes along a subscheme of length $6$ on $S$.}
\end{remark}

\begin{remark}\label{zsection}
{\rm Similarly for a Fano threefold section $X$ of $\Sigma$, the set of 
tangent hyperplane sections $H_t$ that contain $X$ form the $Sp(3)$-dual quartic 
curve $F_{X}$, hence the same curve $F_{X}$ parameterizes rank $2$ vector 
bundles $E$ on $X$ with determinant $H$ and $H^0(X,E)=6$. The general 
section of a vector bundle $E$ vanishes along an elliptic curve of degree 
$6$, a codimension $2$ linear section of a Segre threefold.}
\end{remark}


\begin{lemma}\label{bisecant}
Consider 
a Fano threefold section $X=\Pn {10}\cap\Sigma\subset \PW$. Let 
$C$ be a twisted cubic curve contained in $X$, and let $H_{t}$ be a singular 
hyperplane section of $\Sigma$ that contains $X$.  Then there is a unique 
twisted cubic curve $C^{\prime}\subset X$ 
with length$(C\cap C^{\prime})=2$ determined by $t$.  
In particular, $X$ admits a rank $2$ 
vector bundle $E_{t}$ and a unique section  $s\in \Gamma(X,E_{t})$, 
that $Z(s)=C\cup C^{\prime} $.
\end{lemma}
\begin{proof} Let $Y_{C,t}$ be the unique Segre $3$-fold through $C$ 
that is contained in the tangent hyperplane section $H_{t}$ of 
$\Sigma$.  Then  $Y_{C,t}\cap X=C\cup C^{\prime}$ is a codimension $2$ 
linear section of a Segre $3$-fold.  In particular $C^{\prime}$ is a 
twisted cubic curve and length$(C\cap 
C^{\prime})=2$.  The final statement follows from Proposition \ref{vb} 
and Remark \ref{zsection}.
\end{proof}

Consider the incidence
$$I_{S,T}=\{(\xi,t)\in {\mathcal Hilb}_3S\times T|H^0(S,E_t\otimes I_\xi)>0\}$$
Clearly, by the above 
$$I_{S,T}=\{(\xi,t)|C_\xi\subset H_t\}$$
where $C_\xi$ is the unique twisted cubic curve through $\xi$. 
Now, $T=F_{S}$ is the $Sp(3)$-dual quartic surface to $S$, hence 
$$T_\xi=\{t\in T|(\xi,t)\in I_{S,T}\}$$
is the plane quartic curve on $T$ which is $Sp(3)$-dual 
to the threefold $<S\cup C_{\xi}>\cap \Sigma$.  Therefore 
$$\xi\mapsto T_\xi$$
defines again the map 
$$f:Hilb_3S\to |T_\xi|=\Pn 3_{S}.$$

Let $A_{\xi}$ be the fiber $f^{-1}(f(\xi))$.   Let $\eta\in A_{\xi}$ 
be a general point in the fiber.  By Lemma \ref{bisecant} 
there is a morphism $$\tau_{\eta}: T_{\xi}\to A_{\xi}$$ which assigns to a point 
$t\in T_{\xi}$ the element $(C^{\prime}\cap S)\in A_{\xi}$, where $C^{\prime}$ 
is the unique twisted cubic curve residual to $C_{\eta}$ in the Segre $3$-fold 
defined by $C_{\eta}$ and $t$.  By uniqueness, this map is injective: The 
curve $C_{\eta}$ is the unique twisted cubic curve residual to $C^{\prime}$ 
determined by $t$.  Thus

\begin{proposition}\label{ajm}
    Let $S=\Pn 9\cap \Sigma$ be a smooth linear section with no 
    rational curves of degree less than four. 
    Let $T=F_{S}$ be the $Sp(3)$-dual quartic surface, and consider the 
    abelian fibration $f:Hilb_3S\to |T_\xi|=\Pn 3_{S},$ where 
    $T_{\xi}$ is the plane section of $T$ determined by $C_{\xi}$, 
    the unique twisted cubic curve in $\Sigma$ through $\xi$. 
    Let $A_{\xi}$ be the fiber  $f^{-1}(f(\xi))$. When $T_{\xi}$ is 
    smooth, then any element $\eta\in A_{\xi}$ defines an embedding 
    $\tau_{\eta}: T_{\xi}\to A_{\xi}$.
\end{proposition}

 \section{Twisted cubic curves
          on the prime Fano threefolds of genus $9$}

In the previous section we concluded that the fibers of the fibration 
$$f: Hilb_3S\to \Pn 3_{S}$$
 are abelian threefolds.   From the construction we note that the 
fiber is identified with the subset of ${\mathcal C}(S)$ of twisted cubic curves $C$ 
that are contained in a fixed $\Pn {10}_{h}\supset S$.  In 
particular, the general fiber coincides with the Hilbert scheme of twisted cubic curves in 
the Fano threefold $X_{h}=\Sigma\cap \Pn {10}_{h}$.  The $Sp(3)$-dual 
variety to $X_{h}$ is a plane quartic curve 
$F_{X_{h}}$, the plane section $T_{h}=T\cap\Pn 2_{h}$ of the 
$Sp(3)$-dual surface $T$ to $S$, where $\Pn 2_{h}=(\Pn {10}_{h})^{\bot} \subset\PWD$.
In this section we shall prove the following theorem:

\begin{theorem}\label{abjac}
    Let $S=\Pn 9\cap \Sigma\subset\PW$ be a smooth linear section with no 
    rational curves of degree less than four. 
    Let $T=F_{S}$ be the $Sp(3)$-dual quartic surface, and consider the 
    abelian fibration $f:Hilb_3S\to |T_h|=\Pn 3_{S}.$
The general fiber 
$A_h = f^{-1}(h)$ is isomorphic to the 
jacobian of the $Sp(3)$-dual plane quartic curve 
$T_h$ to the Fano $3$-fold $X_{h}\supset S$.
\end{theorem}

%

For the proof we shall need to know some additional properties 
of the Hilbert scheme ${\mathcal C}(X) $ of twisted cubic curves on the general 
prime Fano threefold $X$ of genus $9$. We begin with:

\begin{proposition}\label{jac}
Let $X =\Pn {10}\cap \Sigma$ be a general prime Fano threefold
of index $1$ and genus $9$, and let $F_X$ be its 
$Sp(3)$-dual plane quartic curve.
Then the intermediate jacobian of $X$ is isomorphic 
to the jacobian $J(F_X)$ of $F_X$. 
\end{proposition}

\begin{proof}
In \cite{Mu3}, Mukai identifies the intermediate jacobian of $X$ 
with the Jacobian of a curve of genus $3$, and in  \cite{Iliev} 
this curve is identified with $F_{X}$.
\end{proof}


We shall show that the family ${\mathcal C}(X)$ is isomorphic to the intermediate 
jacobian $J(X)$. 
For this we shall use the birational properties of the Fano $3$-fold $X$
of genus $9$ related to twisted cubic curves on $X$. 


Let $B \subset X \subset {\bf P}^{10}$ be a smooth rational normal cubic curve, 
and consider the rational projection $$\pi : X \cdots > {\bf P}^6$$ 
from the space ${\bf P}^3 = <B>$.
Let $X' \subset {\bf P}^6$ be the proper $\pi$-image of $X$,
and let $\beta: \tilde{X} \rightarrow X$ be the blowup of $X$ at $B$.
Since we may assume that $<B>\cap X=B$, the blowup $\beta$  
resolves the indeterminacy locus of $\pi$, so the projection $\pi$ 
extends to a morphism 
$$
\tilde{\pi}: \tilde{X} \rightarrow X'.
$$

Recall from \S 4.1 in \cite{IP} that the cubic $B \subset X$ 
fulfills the conditions (*)-(**) if:

\

(*). the anticanonical divisor $-K_{\tilde{X}}$ is numerically effective (nef) 
and $(-K_{\tilde{X}})^3 > 0$ (i.e. $-K_{\tilde{X}}$ is big);

(**) there are no effective divisors $D$ on $X$ such that $(-K_{\tilde{X}})^2.D = 0$. 

\

By the Remark on page 66 in \cite{IP}, the twisted cubic 
$B \subset X$ will satisfy the conditions (*)-(**) if the morphism
$\tilde{\pi}$ has only a finite number of fibers of positive dimension. 

\

\begin{lemma}
Suppose that $X = X_{16}$ is general and the smooth twisted 
cubic $B \subset X$ is general. Then the morphism $\tilde{\pi}$ has only 
a finite number of irreducible fibers of positive dimension, 
and these are precisely the proper 
$\beta$-preimages of the 
$e = e(B) = 12$ 
lines on $X$ that intersect 
the twisted cubic curve $B$.
\end{lemma}

\begin{proof}
First, we may assume that $<B>\cap X=B$ by Lemma \ref{lgsecants} 
and that $B$ is irreducible.  The family of planes in $\Pn 5$ parameterized 
by $B$ sweeps out an irreducible $3$-fold of degree $3$.  
If this $3$-fold is singular, it is a cone, 
and $B\subset Q_{p}$ for some $p\in \Pn 5$, contrary to the assumption.  
Therefore this $3$-fold is smooth, isomorphic to the smooth Segre embedding of 
$\Pn 1\times\Pn 2$. 
Thus any two planes representing points of $B$ are disjoint; 
in particular no conic in $X$ intersects $B$ in a scheme of length $2$, 
since the Lagrangian planes representing the points of a conic in $X$ 
always have a common point.

A positive dimensional fiber of $\tilde{\pi}$ is mapped by $\beta$ 
to a $\Pn 4 \supset <B>$ such that $\Pn 4 \cap X$  
contains a curve $D$ which strictly contains $B$. 

 For a $\Pn 4$ as above the intersection $<D>\cap X$ can't be of dimension $\ge 2$.  
 Indeed the coincidence $<B>\cap X=B$ implies that if $<D>\cap X$ contains a surface $S$, 
 then $S$ would have degree $3$ since its intersection with $<B>$ is $B$.  
But then the span of $S$ must lie in $\Omega$ and $<B>\cap X$ 
must be a quadric surface,
which contradicts to $Pic(X) = {\bf Z}H$.   

So we may assume that $<D>\cap X=D$.  First note that by the ramification 
argument in the proof of Proposition \ref{iso}, 
$D$ is reduced along $B$. Thus $D=B\cup C$, and deg$C=$length$C\cap B$.
  
If deg$C>3$, then $<C>=\Pn 4$, and for $p\in B \setminus C$ 
the projection of $C$ from $p$ will have a $1$-parameter family of trisecant lines. 
Therefore $D = B \cup C$ will have at least a $1$-parameter family of $4$-secant 
planes, and by Lemma \ref{planes} each of these planes will intersect $X$ 
in a conic.  Hence the intersection $X \cap \Pn 4$ will contain at least a surface,
which again contradicts to $Pic(X) = {\bf Z}H$. 

Therefore deg$C \le 3$. 
By the preceding no component of $C$ can be a conic 
(a 1-cycle of degree $2$),  otherwise $X$ will 
have a 2-secant conic which is impossible -- see above. 
Therefore either $C$ is an irreducible twisted cubic or 
any connected component of $C$ is a line. 

But if $C$ is an irreducible irreducible twisted cubic  
then by the preceding the zero-scheme $C\cap B$ must be 
of length 3, i.e. $C$ and $B$ will be two twisted cubic curves 
that pass through the length 3 zero-scheme $C \cap B$, 
thus contradicting to Proposition \ref{iso}. 
 
Therefore any connected component of $C$ must be a line
which is simply secant to $B$, 
i.e. any curve on $X$ contracted by $\pi$ 
must be a secant line to $B$.

Finally, since the Fano 3-fold $X = X_{16}$ is general,
then by Theorem 4.2.7 in \cite{IP}, the family
of lines on $X$ is smooth, 1-dimensional, and the lines
on $X$ sweep out on $X$ a divisor $R \in |{\cal O}_X(4)|$.

Consider a general linear section $S=\Sigma\cap\Pn 9$ 
and the $3$-dimensional family of Fano 3-folds 
$X_t=\Sigma\cap \Pn {10}_{t}$ that contain $S$.
Let $R_t$ be the divisor swept out by the lines on $X_t$, and let  $D_{t}=R_t\cap S$.  
Assume that every twisted cubic on $X_t$ is contained in $R_t$.   
By Theorem \ref{Af} the family of twisted cubic curves on any $X_t$ must be $3$-dimensional, 
by the preceding any of these twisted cubic curves intersects the curve $D_{t}$ 
in a scheme of length $3$, and by Proposition \ref{iso} any length $3$ 
subscheme of $D_{t}$ 
is the intersection on $S$ of a twisted cubic curve on $X_t$.  
Now,  $D_t$ is a curve in
$ |{\cal O}_S(4)|$, and two curves $D_t$ and $D_{t^{\prime}}$ either have
 a common component or intersect each other in a zero scheme of length $256$.  
 Let $Z$ be a subscheme of length $3$ of this intersection.  Then there 
 is a unique twisted cubic curve on both $X_t$ and $X_{t^{\prime}}$, 
 that intersects $S$ in $Z$. But these two curves must be distinct, 
 contradicting to Proposition \ref{iso}.  
 
 Therefore the general twisted cubic $C$ in the general $X = X_t$ lies outside 
 $R = R_t$, in particular $C$ intersects a finite number of lines in $X$.
 Since  $R \in |{\cal O}_X(4)|$ this finite number is $12$.

\end{proof}

\

Let $B \subset X = X_{16}$ be a general twisted cubic curve 
on a general $X = X_{16}$ that satisfies the conclusion of the proposition, 
in particular properties 
(*)-(**). 
Then by Lemma 4.1.1 in \cite{IP}, there exists a commutative diagram

\begin{picture}(100,50)

\put(40,10){\makebox(0,0){$X$}}
\put(80,10){\makebox(0,0){$Y$}}
\put(60,20){\makebox(0,0){$X'$}}
\put(40,40){\makebox(0,0){$\tilde{X}$}}
\put(80,40){\makebox(0,0){$\tilde{X}^+$}}

\put(45,10){\vector(1,0){30}}
\put(40,35){\vector(0,-1){20}}
\put(80,35){\vector(0,-1){20}}
\put(45,35){\vector(1,-1){10}}
\put(75,35){\vector(-1,-1){10}}
\put(45,40){\vector(1,0){30}}

\put(60,07){\makebox(0,0){$\psi$}}
\put(60,43){\makebox(0,0){$\xi$}}
\put(37,25){\makebox(0,0){$\beta$}}
\put(83,25){\makebox(0,0){$\varphi$}}
\put(47,27){\makebox(0,0){$\tilde{\pi}$}}
\put(73,27){\makebox(0,0){$\pi^+$}}
\end{picture}

\noindent
where $\tilde{\pi}$ is as above and  $\pi^+$ is another small contraction,
$\xi: \tilde{X} \cdots > \tilde{X}^+$ 
is a flop 
and $\varphi :\tilde{X}^+ \rightarrow Y$ is a Mori extremal morphism
on the smooth threefold $\tilde{X}^+$, ibid.

\begin{lemma}\label{XY}
(see \cite{IP}Iskovskikh-Prokhorov, Chapter 4, Proposition 4.6.3 (iv)]). 
Let $B \subset X = X_{16}$ be a smooth twisted cubic curve in a smooth $3$-fold which 
fulfills the properties (*)-(**). 
Then in the above diagram the extremal morphism $\varphi: \tilde{X}^+ \rightarrow Y$ 
is of type {\bf E1}, i.e. a contraction of a divisor $D \subset \tilde{X}^+$ 
to a smooth curve $\Gamma \subset Y$ such that: 

\smallskip 

(a) $Y = Y_5 \subset {\bf P}^6$ is a smooth Fano threefold 
of index $r = r(Y) = 2$ and of degree $5$;

\smallskip

(b) $\Gamma = \Gamma^3_9 \subset Y$ is a smooth curve 
of genus $3$ and of degree $9$.  
\end{lemma}

Following the notation from \S 4 in \cite{IP},
we denote by $H$ and $L$ the hyperplane classes of $X$ and $Y$, 
and by $H^* = \beta^*(H)$ and $L^* = \varphi^*(L)$ 
their preimages -- correspondingly in $Pic(\tilde{X})$ 
and in $Pic(\tilde{X}^+)$.  

Let $E = \beta^{-1}(B) \subset \tilde{X}$ be the exceptional 
divisor of $\beta$, and let $E^+ = \xi(E)$ be its proper transform 
in $\tilde{X}^+$. 
Let also $D = \varphi^{-1}(\Gamma) \subset \tilde{X}^+$ 
be the exceptional divisor of $\varphi$, and denote by 
$\tilde{D} \subset \tilde{X}$ its proper transform on $\tilde{X}$,
ibid.: 

\begin{picture}(100,60)

\put(50,10){\makebox(0,0){$X$}}
\put(90,10){\makebox(0,0){$Y$}}
\put(70,20){\makebox(0,0){$X'$}}
\put(50,40){\makebox(0,0){$\tilde{X}$}}
\put(90,40){\makebox(0,0){$\tilde{X}^+$}}
    \put(38,10){\makebox(0,0){$B$}}
    \put(102,10){\makebox(0,0){$\Gamma$}}
    \put(38,40){\makebox(0,0){$E$}}
    \put(102,40){\makebox(0,0){$D$}}
    \put(50,50){\makebox(0,0){$\tilde{D}$}}
    \put(90,50){\makebox(0,0){$E^+$}}

\put(55,10){\vector(1,0){30}}
\put(50,35){\vector(0,-1){20}}
\put(90,35){\vector(0,-1){20}}
\put(55,35){\vector(1,-1){10}}
\put(85,35){\vector(-1,-1){10}}
\put(55,40){\vector(1,0){30}}
    \put(38,35){\vector(0,-1){20}}
    \put(102,35){\vector(0,-1){20}}

\put(70,07){\makebox(0,0){$\psi$}}
\put(70,43){\makebox(0,0){$\xi$}}
\put(47,25){\makebox(0,0){$\beta$}}
\put(93,25){\makebox(0,0){$\varphi$}}
\put(57,27){\makebox(0,0){$\pi'$}}
\put(83,27){\makebox(0,0){$\pi^+$}}
    \put(44,10){\makebox(0,0){$\subset$}}
    \put(96,10){\makebox(0,0){$\supset$}}
    \put(44,40){\makebox(0,0){$\subset$}}
    \put(96,40){\makebox(0,0){$\supset$}}
    \put(50,45){\makebox(0,0){$\cap$}}
    \put(90,45){\makebox(0,0){$\cap$}}

\end{picture}

\begin{proposition}\label{bir}
Let $X = X_{16}$ be general and let $B \subset X$ be a general twisted 
cubic on $X$. Let $\psi: X \cdots > Y$ be the birationality
in the above diagram. 
Denote by ${\mathcal C}(X)$ and by ${\mathcal C}(Y)$ the
families of twisted cubic curves on $X$ and $Y$. 
Then, in the above notation:  

\smallskip

(1.a)  The birational map $\psi: X \rightarrow Y$
is defined by the non-complete linear system $|2H-3B|$. 

\smallskip

(1.b) The extremal divisor of $\psi$, 
the proper transform $D_X \subset X$ of the 
exceptional divisor $D$ of $\varphi$,   
is the unique effective divisor 
of the non-complete linear system $|3H-5B|$ on $X$.  

\smallskip

(1.c) The general extremal curve of $\psi$, the general 
smooth irreducible curve on $X$ contracted by $\psi$, 
is the general element of the compactified family 
${\mathcal C}[2](X;B) \subset {\mathcal C}(X)$ 
of twisted cubic curves $C \subset X$ bisecant to $B$.   

\smallskip

(2.a) The inverse birational map $\psi^{-1}: Y \rightarrow X$ 
is defined by the non-complete linear system $|5L-3\Gamma|$. 

\smallskip

(2.b) The extremal divisor of $\psi^{-1}$,
the proper transform $E_Y \subset Y$ of the 
exceptional divisor $E$ of $\beta$,  
is the unique effective divisor $N$ of the non-complete 
linear system $|3L-2\Gamma|$ on $Y$.  

\smallskip

(2.c) The general extremal curve of $\psi^{-1}$,
the general smooth irreducible curve on $Y$ 
contracted by $\psi^{-1}$, 
is the general element of the compactified family 
${\mathcal C}[5](Y;\Gamma) \subset {\mathcal C}(Y)$ 
of the twisted cubic curves $C \subset Y$ 
that are $5$-secant to $\Gamma$.  
\end{proposition}

\begin{proof} We shall prove (1.a)-(1.b)-(1.c); 
as seen below, the proof of (2.a)-(2.b)-(2.c) is similar. 

\smallskip

{\bf Proof of (1.a)}. 
$\beta$, 
Since $\varphi$ is an extremal morphism of type {\bf E1}, 
then by Lemma 4.1.2 (i) of \cite{IP} and the data 
$g = g(X) = 9$, $g(B) = 0$ and $deg(B) = 3$, one obtains:
$$(-K_{\tilde{X}^+})^3 = 2(g(X)+g(B)-deg(B)-1) - 2 = 8$$
$$(-K_{\tilde{X}^+})^2.E^+ = deg(B) + 2 - 2g(B) = 5$$
$$(-K_{\tilde{X}^+}).(E^+)^2 = 2g(B)-2 = -2.$$
$$E^3=-deg(B)+2-2g(B)=-1$$

%


Since, by Lemma \ref{XY}, the morphism $\varphi: \tilde{X}^+ \rightarrow Y$ is extremal  
of type {\bf E1} onto the Fano threefold $Y$ of index $r = 2$, 
then by Lemma 4.1.5 (ii) of \cite{IP}, 
there exists an integer $\alpha \ge 1$ such that 
$\alpha + 1 \cong 1$(mod r), i.e. $\alpha$ is odd, and   
$$
D = \alpha(-K_{\tilde{X}^+}) - rE^+ = \alpha(-K_{\tilde{X}^+}) - 2E^+
$$
in $Pic(\tilde{X}^+) = {\bf Z}K_{\tilde{X}^+} + {\bf Z}E^+$. 
Next, by taking in mind that $g(\Gamma) = 3$ (see Lemma \ref{XY}), 
the third identity for $\varphi$ on p. 69 in \cite{IP}
gives 
$$
D^2.(-K_{\tilde{X}^+}) = 2g(\Gamma) - 2 = 4.  
$$

Since by the preceding 
$$
D^2.(-K_{\tilde{X}^+})^2 
= (\alpha(-K_{\tilde{X}^+}) - 2E^+)^2.(-K_{\tilde{X}^+})
$$
$$
= \alpha^2 (-K_{\tilde{X}^+})^3 - 4\alpha (-K_{\tilde{X}^+}).E^+ 
  - 4(-K_{\tilde{X}^+})^2.(E^+)^2
= 8\alpha^2 -20 \alpha - 8,
$$
then the identity $D^2.(-K_{\tilde{X}^+}) = 4$ yields 
$$
2\alpha^2 - 5\alpha - 3 = 0.
$$

Since $2\alpha^2 - 5\alpha - 3 = (2\alpha+1)(\alpha-3)$, 
and since by Lemma 4.1.5 (ii) of \cite{IP} the integer $\alpha$ must be positive, 
then the only possibility left is $\alpha = 3$,
i.e. 
$$
D = 3(-K_{\tilde{X}^+}) - 2E^+;
$$
in particular $\alpha$ is odd, ibid.


Since $Y$ is a Fano threefold of index $r = 2$ then $K_Y = -2L$, 
and since $\varphi: \tilde{X} \rightarrow Y$ 
is a blowup of a curve then 
$$
K_{\tilde{X}^+} = \varphi^*(K_Y) + D = -2L^* + D
$$
in $Pic(\tilde{X}^+) = {\bf Z} L^* + {\bf Z}D$; 
notice that $D$ is the exceptional divisor of the 
blowup $\varphi$ of the curve $\Gamma \subset Y$.  
This, together with the previous formula for $D$, yields 
$K_{\tilde{X}^+} = -2L^* + D = -2L^* + 3(-K_{\tilde{X}^+}) - 2E^+$.
Therefore
$$
L^* = 2(-K_{\tilde{X}^+}) - E^+;
$$ 
and since $\xi: \tilde{X} \rightarrow \tilde{X}^+$
is a flop, then in $Pic(\tilde{X})$ 
$$
\tilde{L} = 2(-K_{\tilde{X}}) - E
$$
for the proper $\xi$-preimage $\tilde{L}$ of $L^*$.  

Now we use that $X$ is a Fano threefold of index $1$,
i.e. $-K_X = H$ where $H$ is the hyperplane class of $X$, 
and that $\beta: \tilde{X} \rightarrow X$ is a blowup
of a curve -- the curve $B \subset X$. Therefore 
$$
K_{\tilde{X}} = {\beta}^*(K_{X}) + E = - H^* + E;
$$
recall that $H^* = \beta^*(H)$ is the preimage of $H$ 
and that $E$ is the exceptional divisor of $\beta$. 
Then by the above formulas
$$
L^* = 2(-K_{\tilde{X}}) - E = 2(H^*-E) - E = 2H^* - 3E, 
$$
i.e. the composition map 
$\varphi \circ \xi : \tilde{X} \rightarrow Y$
is given by the linear system  $L^* = |2H^* - 3E|$;
notice that $L$ is the hyperplane class on $Y$.
Since $E$ is the exceptional divisor of the blowup of $B \subset X$, 
the rational map
$\psi: X \rightarrow Y$ is given by the non-complete 
linear system $|2H - 3B|$. 

\medskip

{\bf Proof of (1.b)}. 
We proceed as above:
For the proper $\xi$-preimage $\tilde{D} \subset \tilde{X}$ 
of the exceptional divisor $D \subset \tilde{X}^+$ one obtains
$$
\tilde{D} = 3(-K_{\tilde{X}}) - 2E = 3(H^*-E) - 2E = 3H^* - 5E,
$$
i.e. the extremal divisor $D_X = \beta(\tilde{D}) \subset X$ 
belongs to the non-complete linear system
$|3H - 5B|$ on $X$.

\medskip

{\bf Proof of (1.c)}. 
The general extremal curve $C$ on $X$ is the same as the 
proper transform of the general fiber of $\varphi$. 
Since the map $\tilde{X} \rightarrow Y$ is given by the linear system $\tilde{L}$ 
then the proper preimage on $\tilde{X}$ of such curve   
is the same as the general irreducible curve 
$\tilde{C} \subset \tilde{X}$ such that $\tilde{C}.\tilde{L} = 0$,
i.e. 
$$
\tilde{C}.(2H^* - 3E) = 0.
$$
Therefore the general extremal curve $C \subset X$ 
is a smooth curve of degree $3n$ such that $deg(C.B) = 2n$ 
for some integer $n = n(C) \ge 1$.  
Since  
the compactified family of extremal curves 
on $X$ is a $1$-dimensional algebraic family of $1$-cycles on 
$X$ 
(with a smooth irreducible base isomorphic to the curve $\Gamma$), 
then the integer $n(C) = n$ does not depend on $C$, 
and it remains to see that $n = 1$. 

The equality $n = 1$ is equivalent to say that 
$deg(C.B) = 2$ for the extremal curve $C \subset X$.
By their definition, 
the extremal curves $C \subset X$ are the proper transforms on 
$X$ of the (extremal) curves $C^+ \subset \tilde{X}$ generating the extremal 
ray $R = {\bf R}[C^+]$ defining the extremal morphism $\varphi$. 
Since $E$ is the exceptional divisor of $\beta$ over $B$  
and $E^+$ is its proper transform on $\tilde{X}^+$,
then $deg(C.B) = C^+.E^+$, and it rests to see that 
$C^+.E^+ = 2$. 

For this recall that $\varphi: \tilde{X}^+ \rightarrow Y$ 
is an extremal contraction of type {\bf E1};
in particular the extremal 
ray $R = {\bf R}[C^+]$ (generated by any of the extremal curves
$C^+ \subset \tilde{X}^+$) has length $1$, i.e.
$$
length(R) = -K_{\tilde{X}^+}.C^+ = 1,
$$
see e.g. Theorem 1.4.3 (Mori-Kollar) in \cite{IP}.
Next we use the formula 
$$
E^+ = - 2K_{\tilde{X}^+} -L^*
$$ 
(see above),
and the fact that the extremal curves $C^+ \subset \tilde{X}^+$, 
the fibers of the ruled surface $\varphi|_D: D \rightarrow \Gamma$, 
are orthogonal to the hyperplane class of $Y$, 
i.e. $C^+.L^* = 0$.
Therefore 
$$
C^+.E^+ = C^+.(- 2K_{\tilde{X}^+} - L^*) = 2C^+.(-K_{\tilde{X}^+}) = 2\cdot{\rm length}(R) = 2.
$$
Hence $n = 1$, i.e. the extremal curves $C \subset X$ of $\psi$ are 
the twisted cubic curves on $X$ bisecant to $B$. 

\medskip

The proof of (2.a)-(2.b)-(2.c) is similar; 
notice that $\beta :\tilde{X} \rightarrow X$ is,
just as $\varphi$, an extremal contraction of type {\bf E1}. 
\end{proof} 





\begin{proposition}\label{tc}
The Hilbert scheme ${\mathcal C}(X)$ 
of twisted cubic curves on the general $X = X_{16}$ 
is isomorphic to the jacobian $J(\Gamma)$. 
\end{proposition}

\begin{proof}
By Corollary \ref{cxabel}, the Hilbert scheme ${\mathcal C}(X)$ 
is a $3$-dimensional abelian variety. 
Therefore to prove the isomorphism between the 
abelian varieties ${\mathcal C}(X)$ and $J(\Gamma)$ it will be enough 
to see that they are birationally equivalent. 
For this we shall construct a natural birationality 
$\phi = \phi_{B} : {\mathcal C}(X) \rightarrow Pic^5(\Gamma)$
defined by the choice of a general twisted cubic $B \subset X$. 


Let the twisted cubic $B \subset X$ be as in Proposition \ref{bir},
let $\psi: X \rightarrow Y$ be the birationality defined 
by $B$.   
Let $C \subset X$ be another general twisted cubic on $X = X_{16}
\subset \Pn {10}$; in particular $C$ is smooth, and $C$ spans,
together with $B$, a 7-space $\Pn 7_C = <C \cup B>\subset
\Pn {10}$. In the space $\Pn {10}=<X>$, the family of
codimension 2 subspaces $\Pn 8$ that contain $\Pn 7_C$ is
parameterized by the projective plane $\widehat{\Pn 2_C(t)}$ =
$\Pn {10}\slash\Pn 7_C$. The general $\Pn 8_t \in
\widehat{\Pn 2_C(t)}$ intersects on $X$ a reduced canonical
1-cycle
$$
C^9_{16,t} = B + C + C_{10,t} \subset X
$$
of degree $16$ and of arithmetical genus $9$.
For the general $t$, the residual component $C_{10,t}$ is a smooth
elliptic curve of degree $10$ on $X$ intersecting any of the curves 
$B$ and $C$ at $5$ points.

In the notation of Proposition \ref{bir}, 
the proper $\psi$-image $C_{5,t} \subset Y$
of $C_{10,t} \subset X$ is a smooth projectively normal elliptic
quintic intersecting
$\Gamma = \Gamma^3_9$ at an effective 0-cycle $D_t$ of degree $5$.
When $\Pn 8_t$ moves in $\widehat{\Pn 2_C(t)}$, the
0-cycles $D_t$ describe a $\Pn 2$-family of effective
divisors of degree $5$ on the genus 3 curve $\Gamma$.
Therefore all divisors $D_t$ belong to the same linear
system on $\Gamma$. Since by the Riemann-Roch theorem the dimension
of a complete linear system of degree 5 on a curve of genus 3
is always 2, the family $\{ D_t \}$ is in fact a complete
linear system on $\Gamma$, and we denote this system by $|D_C|$.

Thus we have defined a rational map:
$$
\phi: {\mathcal C}(X) \rightarrow  Pic^5(\Gamma), \
C \mapsto |D_C|.
$$

Next we shall see that the rational map $\phi$ is a birational.
For this we shall construct the rational inverse to $\phi$.




Let $|D| = \{ D_t : t \in \Pn 2 \} \in Pic^5(\Gamma)$
be general, and let $D_t \in |D|$ be general.
The effective $0$-cycle $D_t$ of degree $5$ on
$\Gamma \subset Y$ spans a $4$-space
$\Pn 4_t \subset \Pn 6 = <Y>$. 
Since $Y = Y_5$ is the smooth del Pezzo threefold 
of degree $5$
(i.e. the smooth codimension 3 transversal 
linear section of $G(2,5)$, see Theorem 3.3.1 in \cite{IP}), 
then the 4-space $\Pn 4_t$ intersects on $Y$ a projectively normal
elliptic quintic curve $C_{5,t}$. 

On the del Pezzo threefold $Y$, there exists a $\Pn 1$-family
(a pencil) of quadratic sections
$\{ S_s, s \in \Pn 1 \}$ containing the 1-cycle
$C^8_{14,t} = \Gamma^3_9 + C_{5,t}$
of degree $14$ and of arithmetical genus $8$.
The base locus of this pencil is a half-canonical 1-cycle
$C^{21}_{20} = \Gamma^3_9 + C_{5,t} + C_{6,t}$,
where $C_{6,t}$ is an elliptic sextic on $Y$
intersecting $\Gamma$ at a 0-cycle of degree 9.
Now it is not hard to see
that the residue 1-cycle
$C_{6,t}$ does not depend on the choice of $t$; therefore
$$
C^{21}_{20} = \Gamma^3_9 + C_{5,t} + C_{6},
$$
where $C_6 = C_{6,t}$ for any $t$.

\medskip

By Proposition \ref{bir}, 
the proper $\psi$-preimage of the curve $C_{6} \subset Y$ 
is a twisted cubic curve $C \subset X$,
while the proper preimage of $C_{5,t} \subset Y$
is an elliptic curve $C_{10,t} \subset X$ of degree 10 
intersecting both $B$ and $C$ at $5$ points. 
Turn back to the construction of the rational map $\phi$, 
we see that such $C_{10,t}$ can't be other
than one of the $\Pn 2$-family of residue elliptic curves
of degree 10 in the construction of $\phi(C)$.
Now the construction of $\phi$ implies that $\phi(C) = |D|$.

\medskip 

The existence of $\phi$ and $\phi^{-1}$ at the general points 
shows that the Hilbert scheme ${\mathcal C}(X)$ of twisted cubic curves on $X$
is birational to $Pic^5({\Gamma}) \cong J({\Gamma})$.
This proves the birationality between the abelian threefolds 
${\mathcal C}(X)$ and $Pic^5(\Gamma) \cong J(\Gamma)$, and the proposition 
follows.
\end{proof}

\
Let $F=F_{X}$ be the $Sp(3)$-dual curve
parameterizing singular hyperplane sections of $\Sigma$ 
that contain the Fano threefold $X$, and let $J(F)$ be the jacobian of $F$. 
The embeddings from  Proposition \ref{ajm} translate into embeddings
$$
\tau_{C}:F\to A={\mathcal C}(X) \quad x\mapsto C(x)
$$
defined for a general twisted cubic curve $C$ on $X$ by associating 
to a point $x\in F$ the unique twisted cubic curve $C(x)\subset X$ 
which is bisecant to $C$ and determined by the singular hyperplane section 
of $\Sigma$ defined by $x$.

\begin{proposition}\label{translate1}
The Hilbert scheme ${\mathcal C}(X)$ is isomorphic to the jacobian 
$J(F)$ of the $Sp(3)$-dual quartic curve $F$ to $X$.
Moreover, for the general $C \in A={\mathcal C}(X)$ 
the set $\tau_{C}(F)=F_{C}\subset {\mathcal C}(X)$   
is a translate of the Abel-Jacobi image of the curve 
$F$ in the abelian threefold  $J(F)\cong {\mathcal C}(X)$.
\end{proposition}

\begin{proof}
By the criterion \cite{M} of Matsusaka, for the abelian $3$-fold
$A$ the intersection of any ample divisor $D \subset A$ with
 an effective 1-cycle $F_{C}\subset A$ is at least $3$, 
and the equality $(D\cdot F_{C}) = 3$ is only possible if $A = J(F_{C})$, 
$D$ is a copy of the theta divisor in $J(F_{C})$, and 
$F_{C} \subset A$ is an Abel-Jacobi translate of 
$F_{C}$ in its jacobian $J(F_C) = A$. 

\medskip

Therefore to prove the proposition it will be enough to find 
an ample divisor $D \subset A={\mathcal C}(X)$ such that 
$D\cdot F_{C} = 3$. 

\medskip

In order to describe such a divisor, 
recall that the Fano threefold $X$ contains a 1-dimensional family 
of lines, see e.g. Proposition 4.2.2 and Theorem 4.4.13  and Theorem 4.2.7
in \cite{IP}. In particular this $1$-dimensional family of lines is 
parameterized by a smooth curve of genus $17$. 
 Fix a general line $\lambda \subset X$, 
and consider the subfamily 
$$
D_{\lambda} 
= \{ C \in A: C \cap \lambda \not= \emptyset \}. 
$$
 We first need to show that $ D_{\lambda}$ is an ample divisor.  
 For this consider a curve $E$ on ${\mathcal C}(X)$ and the corresponding union 
 of twisted cubic curves $S_{E}\subset X$.  Since the family of lines in $X$ is parameterized by a smooth curve of genus $17$, 
 the general twisted cubic curve of $E$ cannot intersect every line on $X$.  Therefore the intersection   $
D_{\lambda}\cdot E$ is finite.  On the other hand $S_{E}$ is a hypersurface section of $X$, so the intersection $
D_{\lambda}\cdot E\geq {\rm length}S_{E}\cap \lambda  >0$.  Hence  $ D_{\lambda}$ is ample.

We proceed to show that 
$$
D_{\lambda}\cdot F_{C} = 3
$$
in $A$. 
Now, $F_{C}\subset {\mathcal C}(X)$ is simply the family 
of twisted cubic curves on $X$ bisecant to $C$.
Therefore one has to see that among these 
there are exactly {\it three} twisted 
cubics that intersect the general line $\lambda \subset X$. 

By (1.c) in Proposition \ref{bir}
the family $F_{C}$ of twisted cubic curves
is just the family of extremal curves 
of the birationality $\psi$ defined by the 
twisted cubic $C \subset X$. 
But the extremal twisted cubic curves 
from the family $F_{C}$
sweep out the extremal divisor $D_X \in |3 - 5C|$,
in particular $D_X$ is a cubic hypersurface section of $X$,
see (1.b) in \ref{bir}.   
As a cubic section of $X$, the surface 
$D_X$ intersects the line $\lambda \subset X$ 
at three points.   
Whence there are three twisted cubic curves bisecant to $C $ and secant to $\lambda$ 
 on $X$. 
\end{proof}

\begin{corollary}\label{abjac-cor} 
Let the K3 surface $S = S_{16}$ be general,
and let $T=F_{S}$ be the $Sp(3)$-dual quartic surface to $S$.
Let $h \in \Pn 3_{S}$ be general and let $\Pn 2_h\subset \PWD$
and $\Pn {10}_h\subset\PW$ be the corresponding orthogonal
linear spaces.
Then the Hilbert scheme ${\mathcal C}(X_h)$ of twisted cubic curves on the Fano
$3$-fold $X_h = \LG 36 \cap \Pn {10}_h$  is isomorphic
to the jacobian $J(T_{h})$, where $T_{h}=T\cap \Pn 2_h$. 
\end{corollary}

\begin{proof} 
%
%
%
%
%
First, note that the general $X_h \supset S$ is a general 
prime Fano threefold of genus $9$.  
The hyperplane section $T_h = T \cap \Pn 2_h \subset T$ is exactly 
the $Sp(3)$-orthogonal curve to $X_h$, 
so the corollary follows from Proposition \ref{translate1} .
\end{proof}

Finally, since the general fiber $A_{h}$ of $f$ is isomorphic to 
${\mathcal C}(X_h)$ we have completed the proof of Theorem \ref{abjac}.

\medskip

\section{The group law on the Hilbert scheme of twisted cubic 
curves in a Fano threefold of genus $9$}

We now fix a Fano threefold $X=\Pn {10}\cap \Sigma$ 
and its $Sp(3)$-dual curve $F=F_{X}$ 
parameterizing singular hyperplane sections of $\Sigma$ that contain $X$,
and consider the Jacobian $J(F)$ which by Proposition \ref{translate1} is identified with the
Hilbert scheme  $A={\mathcal C}(X)$ of twisted cubic curves on $X$.
Of crucial importance for our definition of the group law 
on $A={\mathcal C}(X)$ are the embeddings 
$$\tau_{C}:F\to A={\mathcal C}(X) \quad x\mapsto C(x)$$
defined for a general twisted cubic curve $C$ on $X$, by associating 
to a point $x\in F$ the unique twisted cubic curve $C(x)\subset X$ 
bisecant to $C$ and determined by the singular hyperplane section 
of $\Sigma$ defined by $x$.  We denote by $F_{C}\subset A$ the image $\tau_{C}(F)$.

To study the various curves $F_{C}$ in $A$, we fix, once and for all, 
a point $x_{0}\in F$.
Let 
$$
\alpha: F\to F_{0}\subset A\quad x\mapsto C_{x}
$$ 
be the Abel-Jacobi map with base point $x_{0}$, 
and fix the image  $\alpha(x_{0})=C_{0}\in F_{0}\subset A$ 
as the origin of the group structure on $A$.  
Recall from Proposition \ref{vb} that each $x\in F$ 
determines a rank $2$ vector bundle $E_{x}$ on $X$.
\begin{lemma}\label{vbx} For each twisted cubic curve  $C\subset X$, 
the curve $C\cup C(x)$ is
 the zero-locus of a section of the rank $2$ vector bundle $E_{x}$ on $X$.
 \end{lemma}
 \begin{proof}  This follows immediately from Lemma \ref{bisecant}. \end{proof}
\begin{lemma}\label{translate2} 
The translate $F_{\infty} := F_{C} + C$ is independent of $C$. 
\end{lemma}
\begin{proof}
By Lemma \ref{vbx} the point $C+C(x)=C'+C'(x)$ for any $x\in F$, 
since $C\cup C(x)$ and $C'\cup C'(x)$ are both the zero locus 
of a section of the vector bundle $E_{x}$ on $X$.
In particular $F\to A\quad x\mapsto C+C(x)$ is well defined and independent of $C$.  
The image
$\{ C+C(x)| x\in F\}\subset A$ is clearly the translation $F_{C}+C$.
\end{proof}

By the identification of $A$ with $J(F)=Pic_0(F)$, 
there exists a constant $c_{\infty} \in Pic_1(F)$ such that 
the Abel-Jacobi translate 
$$
F_{\infty} = \{ x-c_{\infty}: x \in F \} \subset Pic_0(F) = A,
$$ 
where $x-c_{\infty}$ denotes the rational equivalence class 
of the divisor.
\begin{corollary}\label{sum=x} 
For any $x \in F$ and any $C\in A$, 
the following identity holds
$$
C(x)+C = x - c_{\infty}\in Pic_{0}F=A. 
$$
\end{corollary}

\

\begin{lemma}\label{commutative} 
For any ordered triple $(x,y,z)$ of points on $F$ 
and any $C \in A$, 
let $C(x) \in F_{C}$ be the point representing $x \in F$.
Next, let $C(xy) \in F_{C(x)}$ be the point on $F_{C(x)}$ representing 
the point $y \in F$, 
and let $C(xyz) \in F_{C(xy)}$ be the point representing $z \in F$. 
Then 
$$
C + C(xyz) =  x - y + z - c_\infty = C + C(zyx),
$$
where the addition is taken in $Pic_0(F)$.  
\end{lemma}

\begin{proof} 
By the definitions of $C(x)$, $C(xy)$ 
and $C(xyz)$, 
and by Corollary \ref{sum=x}
$$
C + C(x) = x - c_\infty, \ C(x) + C(xy) = y - c_\infty, 
\mbox{ and }  C(xy) + C(xyz) = z - c_\infty.
$$
Therefore 
$$ C + C(xyz) = (C + C(x) - c_\infty) - (C(x) + C(xy) - c_\infty) + (C(xy) + C(xyz) - c_\infty) $$
 $$= x-y+z-c_\infty.
$$

Similarly $C + C(zyx) = z-y+x-c_\infty = x-y+z-c_\infty$.
\end{proof}

Now we are ready to describe {\it the group law on the threefold $A$}. 
For this we shall need the following elementary technical

\begin{lemma}\label{x-y+z}
For a plane quartic curve $F$, any element $D \in Pic_1(F)$ 
can be represented as the rational equivalence class of a divisor 
$x-y+z$ for some $x,y,z \in F$. 
\end{lemma}

\begin{proof} 
Follows from the Riemann-Roch theorem.
\end{proof}

\

\centerline{\it The sum of two twisted cubic curves in $X$}

\

Let $A$ 
be the abelian $3$-fold above, identified with the jacobian $J(F)$ of the 
plane quartic curve $F$, see Theorem \ref{abjac}.  

Let $(A,\hat{+};C_{0})$ be the additive group structure on 
$A$, 
with the element $C_{0}$ declared to be the neutral element. 
Below we shall find the curve $C''' \in A$ representing the 
sum $C' \ \hat{+} \ C''$ of two general elements 
$C'$ and $C''$ of $A$.

\medskip

By Lemmas \ref{commutative} and \ref{x-y+z},
there exists a triple of points $x,y,z \in F$
such that 
$C'' = C'(xyz)$ is the result of applying 
the chain of transformations 
$$
C' \rightarrow C'(x) \rightarrow C'(xy) \rightarrow C'(xyz)
$$
as above.  Notice that the result is independent of which representative $x-y+z$ in
 its rational equivalence class of divisors one chooses.
Let $C_{0}(xyz)$ be the result of the same chain of transformations to $C_{0}$: 
$$
C_{0} \rightarrow C_{0}(x) \rightarrow C_{0}(xy) \rightarrow C_{0}(xyz),  
$$
i.e. $C_0(x) \in F_{C_{0}}$ is the point representing the point $x \in F$ 
under the isomorphism $F\to F_{C_{0}}$, 
$C_{0}(xy)  \in F_{C_{0}(x) }$ is the point representing the point $y \in F$ 
under the isomorphism $F\to F_{C_{0}(x) }$,
and 
$C_{0}(xyz) \in F_{C_{0}(xy)}$ is the point representing the point $z \in F$ 
under the isomorphism $F\to F_{C_{0}(xy)}$.

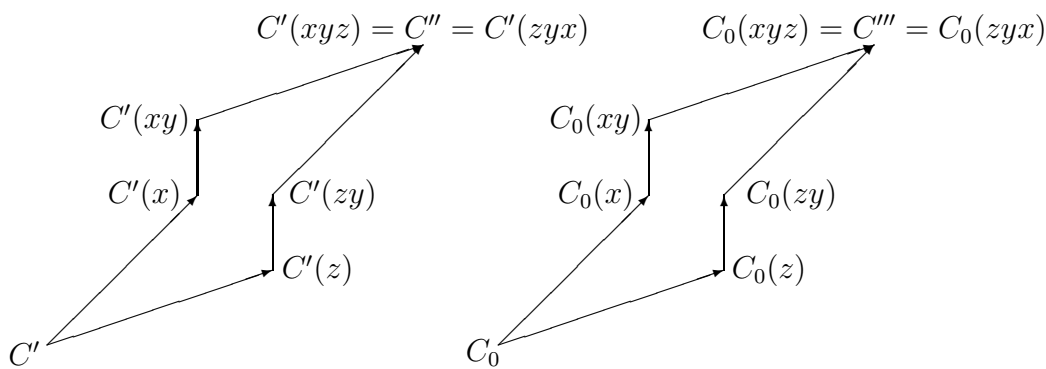
\begin{picture}(120,65)

\put(10,10){\vector(1,1){20}}
\put(10,10){\vector(3,1){30}}
\put(40,20){\vector(0,1){10}}
\put(30,30){\vector(0,1){10}}
\put(40,30){\vector(1,1){20}}
\put(30,40){\vector(3,1){30}}

\put(70,10){\vector(1,1){20}}
\put(70,10){\vector(3,1){30}}
\put(100,20){\vector(0,1){10}}
\put(90,30){\vector(0,1){10}}
\put(100,30){\vector(1,1){20}}
\put(90,40){\vector(3,1){30}}

\put(07,9){\makebox(0,0){$C'$}}
\put(23,30){\makebox(0,0){$C'(x)$}}
\put(23,40){\makebox(0,0){$C'(xy)$}}
\put(46,20){\makebox(0,0){$C'(z)$}}
\put(48,30){\makebox(0,0){$C'(zy)$}}
\put(60,52){\makebox(0,0){$C'(xyz) = C'' = C'(zyx)$}}

\put(68,9){\makebox(0,0){$C_{0}$}}
\put(83,30){\makebox(0,0){$C_{0}(x)$}}
\put(83,40){\makebox(0,0){$C_{0}(xy)$}}
\put(106,20){\makebox(0,0){$C_{0}(z)$}}
\put(109,30){\makebox(0,0){$C_{0}(zy)$}}
\put(120,52){\makebox(0,0){$C_{0}(xyz) = C''' = C_{0}(zyx)$}}

\end{picture}

We shall see that the point 
$$
C''' = C_{0}(xyz) \in A
$$
is the sum of the points $C'$ and $C''$ in $(A_h,\hat{+};C_{0})$. 

For this, we observe that since $C_{0}(xyz)= C'''$ 
is obtained from $C_{0}$ by the same chain of transformations as 
$C'' = C'(xyz)$ from $C'$ then, by Lemma \ref{commutative} 
$$
C' + C'' = C' + C'(xyz)
= x-y+z-c_\infty = C_{0}+ C_{0}(xyz) = C_{0} + C'''.
$$
Therefore 
$C''' + C_{0} = C' + C''$,  
and since $C_{0}$ was chosen to be the {\it zero} in 
$(A_h,\hat{+};C_{0})$
then 
$$
C''' = C' \ \hat{+} \ C'',
$$
in the addition law of the group $(A_h,\hat{+};C_{0})$. 

\

\

\

\

\noindent
Atanas Iliev\\ 
Institute of Mathematics\\
Bulgarian Academy of Sciences\\
Acad. G. Bonchev Str., bl. 8\\
1113 Sofia, Bulgaria\\ 
e-mail: ailiev@math.bas.bg

\

\noindent
Kristian Ranestad\\ 
Matematisk Institutt, UiO\\
P.B. 1053  Blindern\\
N-0316 Oslo, Norway\\ 
e-mail: ranestad@math.uio.no

\end{document}